\documentclass[12pt,leqno,fleqn,epsfig]{article}
\usepackage{amssymb, epsfig, amsmath, amsthm}
\usepackage{mathrsfs}    
\usepackage{color}

\newcommand{\R}{\mathbb{R}}

\newcommand{\blue}{\color{blue}}

\newcommand{\bA}{\boldsymbol A}

\newcommand{\bD}{\boldsymbol D}

\newcommand{\bV}{\boldsymbol V}

\newcommand{\bfsigma}{\boldsymbol\sigma}

\newcommand{\be}{\begin{equation}}
\newcommand{\ee}{\end{equation}}
\newcommand{\bea}{\begin{eqnarray}}
\newcommand{\eea}{\end{eqnarray}}
\newcommand{\bean}{\begin{eqnarray*}}
\newcommand{\eean}{\end{eqnarray*}}

\newcommand{\intl}{\int\limits}

\newcommand{\Beweisende}{\rule{0.2cm}{0.2cm}}

\newcommand{\intmw}{{\int\hspace{-830000sp}-\!\!}}

\newcounter{secnum}

\newtheorem{thm}{Theorem}[section]
\newtheorem{cor}[thm]{Corollary}

\theoremstyle{definition}

\newtheorem{defin}[thm]{Definition}
\newtheorem{rem}[thm]{Remark}

\title{On  Liouville type theorem for a generalized  \\stationary Navier-Stokes equations} 
 
\author{Dongho Chae$^{(*)}$ and J\"{o}rg Wolf$^{(\dagger)}$\\
\ \\
 Department of Mathematics$^{(*), (\dagger)}$ \\
Chung-Ang University\\
Dongjak-gu Heukseok-ro 84\\
Seoul 06974, Republic of Korea\\
and \\
School of Mathematics$^{(*)}$ \\
Korea Institute for Advanced Study\\
Dongdaemun-gu   Hoegi-ro 85 \\
Seoul 02455, Republic of Korea\\
$^{(*)}$e-mail: dchae@cau.ac.kr \\
$^{(\dagger)}$e-mail: jwolf2603@cau.ac.kr}
\date{}
\begin{document}

\maketitle
\begin{abstract}
In this paper we prove a Liouville type theorem for  generalized stationary Navier-Stokes systems in $\Bbb R^3$,    which model  non-Newtonian fluids,  where the  Laplacian term  $\Delta u$ is replaced by 
the corresponding non linear operator $\bA_p( u)=\nabla \cdot ( |\bD(u)|^{p-2} \bD(u))$ with 
$ \bD(u) = \frac{1}{2} (\nabla u + (\nabla u)^{ \top})$, $3/2<p< 3$.  
In the case $3/2< p\le 9/5$ we show that a suitable weak solution $u\in W^{1, p}(\Bbb R^3)$ satisfying  $ \liminf_{R \rightarrow \infty} |u_{ B(R)}| =0$ is trivial, i.e. $u\equiv 0$. On the other hand, for $9/5<p<3$ we impose the condition for the Liouville type theorem in terms of a potential function:  if there exists a matrix valued potential function $\bV$ such that  $ \nabla \cdot \bV =u$, whose $L^{\frac{3p}{2p-3}} $ mean oscillation has  the following growth condition at infinity, $$ 
  \intmw_{B(r)} |\bV- \bV_{ B(r)} |^{\frac{3p}{2p-3}} dx \le C  r^{\frac{9-4p}{2p-3}}\quad  \forall 1< r< +\infty,
 $$
then $u\equiv 0$. In  the  case of  the Navier-Stokes equations, $p=2$, this improves the previous results in the literature.   \ 
  \noindent{\bf
keywords:} stationary Navier-Stokes equations, Liouville type theorem 

\end{abstract}

\section{Introduction}
\label{sec:-1}
\setcounter{secnum}{\value{section} \setcounter{equation}{0}
\renewcommand{\theequation}{\mbox{\arabic{secnum}.\arabic{equation}}}}

We consider the following generalized version of the stationary Navier-Stokes equations in $ \R^{3}$
\begin{align} \label{ns-steady}
 &-\bA_p (u) + (u \cdot \nabla) u = - \nabla \pi\quad  \text{in}\quad  \R^{3},\\
&\qquad\qquad \nabla \cdot u=0,\label{ns-steady1} \end{align} 
 where $u=(u_1, u_2, u_3)=u(x)$ is the velocity field, $\pi =\pi(x)$ is the scalar pressure and 
 $$
\bA_p( u)=\nabla \cdot ( |\bD(u)|^{p-2} \bD(u)), \quad  1 < p< +\infty
  $$
with  $ \bD(u) =\bD= \frac{1}{2} (\nabla u + (\nabla u)^{ \top})$ representing  the symmetric gradient. 
 Here $  |\bD|^{p-2} \bD= \bfsigma (\bD)$ stands for the deviatoric stress tensor.
 The system  \eqref{ns-steady}-\eqref{ns-steady1}  is popular among engineers, known as a power law model of  {\em non-Newtonian fluid}, where the  viscosity depends on the shear rate $ |\bD(u)|$.    
  For $p=2$ it reduces to the usual stationary Navier-Stokes equations. For $ 1< p< 2$ the fluid is called shear thinning, while  in case $ 2< p< +\infty$ the fluid is called shear thickening.  For more details on the continuum mechanical  
  background  of the above equations we refer to  \cite{wil}.  
 
\vspace{0.2cm}

The  Liouville type problem for the Navier-Stokes equations, as stated in  Galdi's book\cite[Remark X. 9.4, pp. 729]{gal}, is a challenging open problem in the  mathematical fluid mechanics. We refer \cite{ser1, ser2, ser3, koc, cha2, cha3, ser3, koz, cham, pil, gil} and  the references therein for partial progresses for the problem.  In those literatures authors provided sufficient conditions for velocities to guarantee the triviality of solutions.\\

We say  $\bV \in L^1_{loc} (\R^{n}; \R^{n\times  n}  )$ is  a potential function for  vector field $u\in L^1_{loc}  (\Bbb R^n)$ if  $ \nabla \cdot \bV=u$, where the derivative is in the sense of distribution. 
In \cite{ser1, ser2} Seregin proved Liouville type theorem for the Navier-Stokes equations  under  hypothesis  on the potential function $\bV$ for a solution $u$. In particular in \cite{ser2} it is shown that if $\bV \in BMO (\Bbb R^3)$, then $u=0$. In this paper we would like to improve and generalize  this result for  the system \eqref{ns-steady}-\eqref{ns-steady1}.\\

 For a measurable set $\Omega \subset \Bbb R^n$  we denote by $|\Omega |$ the $ n$-dimensional Lebesgue measure of $\Omega$, and for $f\in L^1 (\Omega)$ we use the notation
$$ f_\Omega :=\intmw_{\Omega}  f dx:= \frac{1}{|\Omega|} \int_{\Omega} fdx.
$$

In contrast to the case $ p=2$ it is still open whether any weak solution to the system  \eqref{ns-steady}-\eqref{ns-steady1} is regular or not. Therefore, in the present paper we only work with  weak solutions satisfying the local energy inequality the solution of which are called suitable weak solution. 

 \begin{defin}
 \label{def1.1}
Let $ \frac{3}{2} \le  p < +\infty$.  

1. We say $ u\in W^{1,\,p}_{ loc}(\R^{3} ) $ is a weak solution to  \eqref{ns-steady}-\eqref{ns-steady1} if the following identity is fulfilled 
\begin{equation}
 \intl_{ \R^{3} } \Big(|\bD(u)|^{ p-2} \bD(u)  -  u \otimes u\Big) : \bD(\varphi ) dx = 0 
\label{weak}
 \end{equation}  
  for all  vector fields $ \varphi \in C^\infty_c(\R^{3} )$ with 
$ \nabla \cdot \varphi =0$.

2. A pair $ (u, \pi ) \in W^{1,\,p}_{ loc} ( \R^{3} ) \times  L^{ \frac{3}{2}}_{ loc} (\R^{3} )$ is called a suitable weak solution to \eqref{ns-steady}, 
\eqref{ns-steady1} if besides  \eqref{weak} 
the following local energy inequality holds 
\begin{align}
 &\intl_{ \R^{3} } |\bD(u)|^{ p} \phi  dx  
\cr
& \le  \intl_{ \R^{3} } |\bD(u)|^{  p-2} \bD(u): u \otimes \nabla \phi   dx + 
\intl_{ \R^{3} } \Big(\frac{1}{2} |u|^2 + \pi \Big)  u \cdot \nabla \phi    dx
\label{sweak}
\end{align}
for all non-negative $ \phi \in C^\infty_{ c}(\R^{3} )$.

\end{defin}

 \begin{rem}
 \label{rem1.2}  In case $ \frac{9}{5} \le p < +\infty$ any weak solution to \eqref{ns-steady}-\eqref{ns-steady1} is a suitable weak solution. 
 Indeed,  by Sobolev's embedding theorem we have  $ u\in L^{ \frac{9}{2}}(\R^{3} )$,  which yields 
$ |u|^2 |\nabla u| \in L^1_{ loc}(\R^{3} ) $.  In addition, as we will see below in Section\,2  from  \eqref{weak} we get  
$ \pi \in L^{ \frac{9}{4}}_{ loc}(\R^{3} )$  such that for all $ \varphi \in W^{1,\, \frac{9}{5}} (\R^{3})$ with compact support 
\begin{equation}
 \intl_{ \R^{3} } \Big(|\bD(u)|^{  p-2} \bD(u): \bD(\varphi )  +  u \otimes u : \bD(\varphi ) \Big)dx =
  \intl_{ \R^{3} } \pi  \nabla \cdot \varphi dx.   
\label{weak1}
 \end{equation}   
 Thus,  inserting $ \varphi = u\phi $ into  \eqref{weak1}, where $ \phi \in  C^\infty_c( \R^{3} )$, and applying integration by parts, we get 
  \eqref{sweak} where the inequality is replaced by equality.

\end{rem}

Our aim in this paper is to prove the following.

\begin{thm}
 \label{thm1.1}
 \begin{itemize}
 \item[(i)]  \underline{ Let $ \frac{3}{2} \le p \le  \frac{9}{5}$}.  We suppose $ (u, \pi ) \in  W^{1,\,p}_{  loc} (\R^{3} )\times L^{ \frac{3}{2}}_{ loc}(\R^{3} )$ is a suitable weak  solution of  \eqref{ns-steady}-\eqref{ns-steady1}.  If 
\begin{equation}
 \intl_{ \R^{3} } |\nabla u|^p dx <+\infty,\quad  \liminf_{R \rightarrow \infty} |u_{ B(R)}| =0
\label{1.2b}
 \end{equation} 
then $ u \equiv 0$. 
\item[(ii)]  \underline{Let $ \frac{9}{5} < p <3$}. We suppose  $ (u, \pi ) \in  W^{1,\,p}_{ loc} (\R^{3} )\times L^{ \frac{3}{2}}_{ loc}(\R^{3} )$ is  a  weak  solution of  \eqref{ns-steady}-\eqref{ns-steady1}.  
Assume there exists $ \bV\in W^{2,p}_{loc }(\R^{3}; \R^{3\times  3}  )$ such that 
 $ \nabla \cdot \bV = u$, and 
 \begin{equation}
  \intmw_{B(r)} |\bV- \bV_{ B(r)} |^{\frac{3p}{2p-3}}dx \le C r^{\frac{9-4p}{2p-3}}\quad  \forall 1< r< +\infty.
 \label{1.2}
  \end{equation} 
 Then,  $ u \equiv 0 $.
 \end{itemize}
  \end{thm}
\begin{rem}
 \label{rem1.4}
Obviously $ \bV \in BMO  (\R^{3} )$ implies the condition  \eqref{1.2}. In fact,  \eqref{1.2}   
is guaranteed  by $ \bV \in C^{ 0,\alpha}(\R^{3} )$  wih $\alpha=\frac{ 9-4p}{3p} >0$ thanks to the Campanato theorem\cite{gia}.      \\
\ \\
As an  immediate corollary of the above theorem we have the following result, which is  the case of $p=2$, which improves the previous result in \cite{ser1, ser2}.
\begin{cor}
 \label{cor1.1}
 Let $ (u, \pi )$ be a smooth solution of   the  stationary Navier-Stokes equations on $\Bbb R^3$.  Suppose there exists $ \bV\in  C^\infty(\R^{3}; \R^{3\times  3}  )$ such that 
 $ \nabla \cdot \bV = u$, and 
 \begin{equation}
  \intmw_{B(r)} |\bV- \bV_{ B(r)} |^{6}dx \le C r \quad  \forall 1< r< +\infty.
 \label{1.2a}
  \end{equation} 
 Then,  $ u \equiv 0 $.
\end{cor}

\end{rem}  
 \section{Proof  of Theorem\,\ref{thm1.1}}
 \label{sec:-2}
 \setcounter{secnum}{\value{section} \setcounter{equation}{0}
 \renewcommand{\theequation}{\mbox{\arabic{secnum}.\arabic{equation}}}}

 We start our discussion of estimating the  pressure for both of the cases (i) and (ii).  First note that  by the hypothesis  $u\in W^{1,p}_{loc} (\Bbb R^3)$ and due to Sobolev's embedding theorem it holds 
$ u \in L^{ \frac{3p}{3-p}}(\R^{3} )$.   This yields 
\[
|\bD(u)|^{  p-2} \bD(u)  -  u \otimes u\in L^{q}_{ loc} (\R^{3} ), \quad  q = \min\Big\{  \frac{3p}{6-2p}, \frac{p}{p-1}\Big\}.   
\]
Given  $ 0 <  R< +\infty$, and noting that $ q \ge \frac{3}{2}$ for $p\ge \frac{3}{2}$, we may  define the functional   
$ F\in W^{-1,\, s}(B(R))$,  $ \frac{3}{2} \le s  \le q$,  by means of  
\[
 \langle F, \varphi \rangle =  \intl_{B(R)} (|\bD( u)|^{p-2} \bD( u)  - u \otimes u) : \bD(\varphi) dx,\quad  \varphi \in W^{1,\, s'  }_0(B(R)),
\]
 where we set $s'=\frac{s}{s-1}$.
Since $ u$ is a  weak  solution to  \eqref{ns-steady}-\eqref{ns-steady1} in view of  \cite[Lemm{\blue a}\,2.1.1]{sohr}   there exists a unique 
$ \pi_R\in L^{ q}(B(R))$ with  $ \intl_{B(R)} \pi _R dx =0$ such that 
  \[
 \langle F, \varphi \rangle = \intl_{B(R)} \pi_R\nabla \cdot \varphi dx \quad \forall \varphi \in W^{1,  s' }_0(B(R)).
\]
 Furthermore,  we get for all $ \frac{3}{2} \le s \le q$ 
\begin{align}
 \label{2.2a1}\intl_{B(R)} |\pi_R|^{s} dx   &\le c \|F\|^{s}_{ W^{-1,\,s} (B(R))} \le  
 c\||\bD(u)|^{p-2} \bD(u) - u \otimes u\|^s_{ L^{s}(B(R))},
\end{align}
 with a constant $ c>0$, depending only on $ p$ but independent of $ 0< R< +\infty$. Let $ 1 {< \blue  \rho<} R <+\infty$. We  set $ \widetilde{\pi}_R= \pi _R - (\pi _R)_{ B(1)} $. 
From the definition of the  pressure $ \pi _R$ it follows that  
\[
\intl_{B(\rho )} (\widetilde{\pi} _R - \widetilde{\pi} _{ \rho }) \nabla \cdot \varphi dx =0 \quad \forall \varphi \in W^{1,  s'  }_0(B(\rho )).
\]
This shows that $\widetilde{\pi}_R - \widetilde{\pi}_\rho$  is constant in $ B(\rho )$.  Since $ (\widetilde{\pi}_R - \widetilde{\pi}_\rho)_{ B(1)}=0$ it follows that 
$ \widetilde{\pi} _{ \rho } = \widetilde{\pi}_R $ in $ B(\rho )$.  This allows us to define $ \pi \in L^q_{ loc}(\R^{3} )$ by setting 
$ \pi = \widetilde{\pi} _{ R}$ in $ B(R)$.  In particular,  $ \pi-\pi_{ B(R)} = \pi_R$. Thus, thanks to  \eqref{2.2a1} we estimate by  H\"older's inequality 
\begin{align*}
\intl_{B(R)} |\pi - \pi _{ B(R)}|^{s} dx &\le c\| |\bD(u)|^{p-2} \bD(u) - u \otimes u\|^{ s}_{ L^{s}(B(R))}
\\
&\le c  R^{ \frac{s(3-p)}{p}} \bigg(\intl_{B( R)} |\bD( u)|^p dx \bigg)^{ s(p-1)} + c  \intl_{B( R)} |u|^{ 2s} dx.
\end{align*}
Hence,  
\begin{equation}
 \|\pi-\pi_{ B(R)}\|_{L^{ s}(B(R))} \le c R^{ \frac{3-p}{p}} \|\bD(u)\|_{ L^p(B(R))}^{ p-1} + c\|u\|^2_{ L^{ 2s}(B(R))}.  
\label{2.21b}
 \end{equation} 
Note that $ q = \frac{9}{4}$ whenever $ \frac{9}{5} \le p < +\infty$. This yields the existence of the pressure 
$ \pi \in L^{ \frac{9}{4}}_{ loc}(\R^{3} )$.

Let $ 1 < r< +\infty$ be arbitrarily chosen, and  $ r \le  \rho < R \le 2r $. We set $ \overline{R}= \frac{R+ \rho }{2}$.  Let $ \zeta \in  C^\infty(\R^{n} )$ be a cut off function, which is radially non-increasing   with $ \zeta =1$ on $B(\rho )$ and 
$\zeta=0$ on $ \Bbb R^3\setminus B(\overline{R})$ satisfying 
$ |\nabla \zeta | \le c (R-\rho )^{ -1}$. From  \eqref{sweak} with $ \phi = \zeta ^p$ we get 
\begin{align*}
 \intl_{B(R )} |\bD( u)|^p \zeta ^p dx & \le  \intl_{B(\overline{R})}  |\bD( u)|^{p-2} \nabla  \zeta ^p \cdot  \bD( u)  \cdot u dx + 
\\
 &\qquad + \frac{1}{2}\intl_{B(\overline{R})} |u|^p u \cdot \nabla \zeta ^p + 
 \intl_{B(\overline{R})} (\pi - \pi_{ B(\overline{R})}) u \cdot \nabla \zeta ^p dx.  
 \end{align*}
Applying H\"older's and Young's inequality, we get from above 
 \begin{align}
 \intl_{B(\rho )} |\bD( u)|^p \zeta ^p dx & \le    c(R-\rho )^{-p}   \intl_{B(\overline{R})\setminus B(\rho )}  |u|^p dx +  c (R- \rho )^{ -1} \intl_{B(\overline{R}) \setminus B(\rho )} |u|^3  dx 
 \cr
 &\qquad + c   (R- \rho )^{ -1} \intl_{B(\overline{R})\setminus B(\rho )} |\pi  - \pi _{ B(\overline{R}) }| |u| dx\cr
 &=   I+ II+ III. \label{2.1}
\end{align}

 \underline{\it The case $ \frac{3}{2} \le p \le \frac{9}{5}$}:   Observing  \eqref{1.2a} and applying Sobolev's embedding theorem we get 
\begin{equation}
u\in L^{ \frac{3p}{3-p}}(\R^{3} ).  
\label{2.1a}
 \end{equation} 
In  \eqref{2.1} we take  $ \rho = \frac{R}{2}$.  Applying H\"older's inequality, we easily get 
\[
I+ II \le  c  \bigg(\intl_{ \R^{3}  \setminus B( \frac{R}{2})} |u|^{ \frac{3p}{3-p}} dx\bigg)^{ \frac{3-p}{3}} 
+ c R^{ \frac{5p-9}{p}}  \bigg(\intl_{ \R^{3}  \setminus B( \frac{R}{2})} |u|^{ \frac{3p}{3-p}} dx\bigg)^{ \frac{3-p}{p}}. 
\]
Using  \eqref{2.1a} and recalling that $ p \le \frac{9}{5}$, we see that $ I+II= o(R)$ as $ R \rightarrow +\infty$.  

Applying H\"older's inequality along with  \eqref{2.21b} with $ s= \frac{3}{2}$,  we estimate 
\begin{align*}
III &\le cR^{ -1}  \Big(R^{ \frac{3-p}{p}} \|\bD(u)\|_{ L^p(B(\overline{R}))}^{ p-1} + c\|u\|^2_{ L^3(B(\overline{R}))}\Big) 
\bigg(\intl_{ \R^{3}  \setminus B( \frac{R}{2})} |u|^{ 3} dx\bigg)^{ \frac{1}{3}}
\\
& \le c \|\nabla u\|_{ L^p}^{ p-1}    \bigg(\intl_{ \R^{3}  \setminus B( \frac{R}{2})} |u|^{ \frac{3p}{3-p}} dx\bigg)^{ \frac{3-p}{3p}} 
+ c R^{ \frac{5p-9}{p}} \|u \|^2_{ L^{ \frac{3p}{3-p}}} \bigg(\intl_{ \R^{3}  \setminus B( \frac{R}{2})} |u|^{ \frac{3p}{3-p}} dx\bigg)^{ \frac{3-p}{3p}}. 
\end{align*}
Observing   \eqref{2.1a} along with  $ p \le \frac{9}{5}$,  we find  $ III= o(R)$ as $ R \rightarrow +\infty$.  Inserting the above estimates into the right-hand side 
of  \eqref{2.1}, we deduce that $ \bD(u) \equiv 0$, which implies that  $ u=u(x) $ is a linear function $x$.  Taking into account the condition   \eqref{1.2b}, we obtain $ u \equiv 0$.

\vspace{0.3cm}
\underline{\it The case $ \frac{9}{5} < p < 3$}:    In order to estimate $ I$ and $ II$ we choose another   cut off function  $ \psi \in  C^\infty(\R^{3} )$, 
which is radially non-increasing   with $ \psi =1$ on $ B(\overline{R})$ and 
$\psi=0$ on $ \Bbb R^3\setminus B(R)$ satisfying 
  $ |\nabla \psi  | \le c (R-\rho )^{ -1}$.
   Recalling that $ u = \nabla \cdot \bV$, applying integration by parts   and  applying the H\"{o}lder inequality,  
we find 
\begin{align*}
\intl_{B(R)} |u|^p\psi ^p dx &= \intl_{B(R)} \partial _i (V_{ ij}- (V_{ ij})_{ B(R)}) u_{ j}  |u|^{p-2}\psi^p dx  
\\
&  =  -\intl_{B(R)} (V_{ ij}- (V_{ ij})_{ B(R)})\bigg(  \partial _i u_{ j} |u|^{p-2}  +  (p-2) u_j u_k \partial _i u_k |u|^{p-4} \bigg)  \psi^p dx\cr
&\qquad -  \intl_{B(R)} (V_{ ij}- (V_{ ij})_{ B(R)})u_j |u|^{p-2}\partial _i\psi^p dx\cr
& \le c\bigg(\intl_{B(R)}| \bV - \bV_{ B(R)}|^p dx \bigg)^{  \frac{1}{p}} \bigg(\intl_{B(R)}| \nabla u|^p dx \bigg)^{  \frac{1}{p}}   \bigg(\intl_{B(R)}| u|^p\psi ^p dx \bigg)^{  \frac{p-2}{p}}\\
&\qquad + c (R- \rho )^{ -1}\bigg(\intl_{B(R)}| \bV - \bV_{ B(R)}|^p dx \bigg)^{  \frac{1}{p}} \bigg(\intl_{B(R)}| u|^p\psi ^p dx \bigg)^{  \frac{p-1}{p}}.
\end{align*}   
 Using H\"{o}lder's inequality, Young's inequality, and  observing  \eqref{1.2}, we obtain 
 \begin{align*}
   \intl_{B(R)} |u|^p\psi ^p dx & \le c \bigg(\intl_{B(R)}| \bV - \bV_{ B(R)}|^p dx \bigg)^{  \frac{1}{2}} \bigg(\intl_{B(R)}| \nabla u|^p dx \bigg)^{  \frac{1}{2}}
\\
&\qquad + c   (R- \rho )^{ -p}\intl_{B(R)}| \bV - \bV_{ B(R)}|^p dx \cr
\end{align*}
\begin{align*}
&\le    cR^{3-p} \bigg(\intl_{B(R)}| \bV - \bV_{ B(R)}|^{\frac{3p}{2p-3}}  dx \bigg)^{  \frac{2p-3}{6}} \bigg(\intl_{B(R)}| \nabla u|^p dx \bigg)^{  \frac{1}{2}}   \cr
&\qquad + c   (R- \rho )^{ -p} R^{6-2p} \bigg(\intl_{B(R)}| \bV - \bV_{ B(R)}|^{\frac{3p}{2p-3}}  dx \bigg)^{  \frac{2p-3}{3}} \cr
&\le c R^{\frac{9-2p}{3}} \bigg(\intl_{B(R)}| \nabla u|^p dx \bigg)^{  \frac{1}{2}}  + c (R-\rho)^{-p} R^{\frac{18-4p}{3}}.
 \end{align*}  
Since $ R\ge1$, and $ p>9/5$ we have $ R^{\frac{9-2p}{3}} \le  R ^p$ and  $ R^{\frac{18-4p}{3}} \le R^{2p}$, and therefore
 \begin{align*}
 I &\le  c  (R- \rho )^{ -p}  R^p  \bigg(\intl_{B(R)}| \nabla u|^p dx \bigg)^{  \frac{1}{2}} +  (R-\rho )^{ -2p} R^{2p}.
 \end{align*}

To estimate $ II$ we proceed similar. We first estimate the $ L^3$ norm of $ u$ as follows 
\begin{align*}
\intl_{B(R)} |u|^3 \psi^3 dx &= \intl_{B(R)} \partial _i (V_{ ij}- (V_{ ij})_{ B(R)}) u_{ j} |u|\psi^3 dx  
\\
&  = - \intl_{B(R)} (V_{ ij}- (V_{ ij})_{ B(R)}) \partial _i (u_{ j} |u|)\psi^3 dx -    \intl_{B(R)} (V_{ ij}- (V_{ ij})_{ B(R)}) u_{ j}|u| \partial _i\psi^3 dx
\\[0.2cm]
& \le c\bigg(\intl_{B(R)}| \bV - \bV_{ B(R)}|^{\frac{3p}{2p-3}}  dx \bigg)^{  \frac{2p-3}{3p}} \bigg(\intl_{B(R)}| u|^3 \psi ^3 dx \bigg)^{  \frac{1}{3}}\bigg(\intl_{B(R)}| \nabla u|^p dx \bigg)^{  \frac{1}{p}}
\\
&\qquad + c (R- \rho )^{ -1}\bigg(\intl_{B(R)}| \bV - \bV_{ B(R)}|^3 dx \bigg)^{  \frac{1}{3}} \bigg(\intl_{B(R)}| u|^3\psi ^3 dx \bigg)^{  \frac{2}{3}}.
\end{align*}   
Using  Young's inequality, we get 
\begin{align}
\intl_{B(R)} |u|^3\psi ^3 dx 
& \le c\bigg(\intl_{B(R)}| \bV - \bV_{ B(R)}|^{\frac{3p}{2p-3}} dx \bigg)^{  \frac{2p-3}{2p}}\bigg(\intl_{B(R)}| \nabla u|^p dx \bigg)^{  \frac{3}{2p}}\cr
& \qquad + c (R- \rho )^{ -3} \intl_{B(R)}| \bV - \bV_{ B(R)}|^3 dx\cr
& \le c\bigg(\intl_{B(R)}| \bV - \bV_{ B(R)}|^{\frac{3p}{2p-3}} dx \bigg)^{  \frac{2p-3}{2p}}\bigg(\intl_{B(R)}| \nabla u|^p dx \bigg)^{  \frac{3}{2p}}
\cr
&\qquad + c (R- \rho )^{ -3}   R^{\frac{3(3-p)}{p}}\bigg(\intl_{B(R)}| \bV - \bV_{ B(R)}|^{\frac{3p}{2p-3}} dx \bigg)^{  \frac{2p-3}{p}} .
 \label{2.3}
\end{align}   
Once more appealing to  \eqref{1.2}, and recalling $ R \ge 1$, $p>9/5$,  and  thus $R^{\frac{9-p}{p}}  \le R^4$,  we arrive at 
 \begin{align}\label{2.3a}
 II & \le  c (R- \rho )^{ -1} R\bigg(\intl_{B(R)}| \nabla u|^p dx \bigg)^{  \frac{3}{2p}} + c  (R- \rho )^{ -4} R^{\frac{9-p}{p}} \cr
&\le c (R- \rho )^{ -1} R\bigg(\intl_{B(R)}| \nabla u|^p dx \bigg)^{  \frac{3}{2p}} +c  (R- \rho )^{ -4} R^4. \end{align}  
It remains to estimate $ III$. Using H\"older's inequality and Young's inequality,  we infer 
\begin{align}
III &\le c (R- \rho )^{ -1} \intl_{B(\overline{R})} |\pi- \pi_{ B(\overline{R})}|^{ \frac{3}{2}} dx   + 
c (R- \rho )^{ -1} \intl_{B(\overline{R})} |u|^{3} dx.   
 \label{2.20} 
\end{align}

Combining  \eqref{2.20},  \eqref{2.3a} and  \eqref{2.21b},  we obtain 
\[
III \le c   R^{ \frac{3(3-p)}{2p}} (R- \rho )^{ -1} \bigg(\intl_{B( \overline{R})} |\nabla u|^p dx \bigg)^{ \frac{3(p-1)}{2p}} + c  (R- \rho )^{ -1}\intl_{B(\overline{R})} |u|^3 dx.
\]

The second   term on the right-hand side can be absorbed into  $ II$. We also observe here,  $R^{ \frac{3(3-p)}{2p}} <R$ thanks to $R\ge 1$ and $p>9/5$.

Thus, inserting the estimate of $ II$,   and once more using $ R \ge 1$,   we find 
 \begin{align*}
 III &\le   c   R (R- \rho )^{ -1}  \bigg(\intl_{B( \overline{R})} |\nabla u|^p dx \bigg)^{ \frac{3(p-1)}{2p}} +     c R(R- \rho )^{ -1}\bigg(\intl_{B(R)}| \nabla u|^p dx \bigg)^{  \frac{3}{2p}} \cr
 & \qquad +  cR^{4} (R-\rho )^{ -4}.
  \end{align*}  
Inserting the estimates of $ I, II$ and $ III$ into the right hand side of  \eqref{2.1}, and applying Young's inequality,  we are led to  
  \begin{align}\label{iter}
\intl_{B(R)} | \bD( u)|^p \zeta^p dx    &\le \frac12 \intl_{B(R )} |\nabla u|^p  dx 
 +  cR^{4} (R-\rho )^{ -4}+      c R^{2p} (R-\rho)^{-2p} \cr
 & \qquad+c R^{ \frac{ 2p}{2p-3}} (R-\rho)^{ -\frac{ 2p}{2p-3}}  +c R ^{\frac{2p}{3-p} } (R-\rho) ^{-\frac{2p}{3-p} } \cr
 &\le \frac12 \intl_{B(R )} |\nabla u|^p dx 
 + c R^m(R-\rho)^{-m},   \end{align}
 where we set
 $$ m= \max\bigg\{ 4, 2p, \frac{ 2p}{2p-3},  \frac{2p}{3-p}  \bigg\},
 $$
{\color{red} and }   used the fact that  $ R^\alpha  (R-\rho)^{-\alpha} \le R^\beta (R-\rho )^{-\beta} $ for $\alpha \le \beta $.
  Furthermore, applying Calder\'on-Zygmund's inequality, we infer 
 \begin{align}
\intl_{B(\rho )} |\nabla u|^p  dx &\le  \intl_{ \R^{3} } |\nabla (u \zeta) |^p  dx    
 \cr
 &\le  \intl_{B(R )} | \bD (u)  |^p \zeta ^p  dx    
+ c(R-\rho )^{-p }  \intl_{B(R)} | u |^p  dx.  
 \label{korn} 
 \end{align}
 Estimating the left-hand side of  \eqref{iter} from below by  \eqref{korn},  
and applying the  iteration Lemma in \cite[V.\,Lemma\,3.1]{gia},  we deduce that 
\begin{align}\label{2.4a}
\intl_{B(\rho ) } |\nabla u|^p dx &\le c R^m (R-\rho)^{-m}
\end{align} 
for all $ r\le \rho <R \le 2r$.
Choosing $ R=2r$ and $\rho=r$ in \eqref{2.4a}, and passing $r\to +\infty$, we find 
\begin{equation}
\int_{\Bbb R^3} |\nabla u|^p dx<+\infty .
\label{2.4}
 \end{equation} 
 Similarly, from  \eqref{2.3a}  and \eqref{2.4}, we get the estimate 
 \begin{equation}
 r^{ -1}\intl_{B(r)} |u|^3 dx \le c  \quad  \forall 1< r< +\infty. 
 \label{2.5}
  \end{equation}

Next, we claim that 
\begin{equation}
 r^{ -1} \intl_{B(3r) \setminus B(2r)} |u|^3 dx = o(1)\quad  \text{as}\quad r \rightarrow +\infty. 
\label{2.6}
 \end{equation} 
 Let $ \psi \in  C^\infty( \R^{3} )$ be a cut off function for the annulus $ B(3r) \setminus B(2r)$ in $ B(4r) \setminus B(r)$, 
i.e. $ 0 \le \psi \le 1$ in $ \R^{3} $, $ \psi = 0 $ in $ \R^{3} \setminus (B(4r) \setminus B(r))$, $ \psi =1$ on $ B(3r) \setminus B(2r)$ and  
 $ |\nabla \psi | \le c r^{-1}$.  Recalling   that $ u = \nabla \cdot \bV$,    and applying integration by parts, using H\"older's inequality along with  \eqref{1.2} 
 we calculate 
 \begin{align*}
&  \intl_{B(4r) \setminus B(r)} |u|^3 \psi^3   dx   
\cr
& \quad = \intl_{B(4r) \setminus B(r)} \partial _j (V_{ ij}- (V_{ ij})_{ B(4r)}) u _i |u| \psi^3   dx  
\cr
&\quad  = -\intl_{B(4r) \setminus B(r)}  (V_{ ij} - (V_{ ij})_{ B(4r)}) \partial _j (u _i |u|) \psi^3   dx  -  
\intl_{B(4r) \setminus B(r)}   (V_{ ij}- (V_{ ij})_{ B(4r)}) (u _i |u|)\partial _j \psi^3   dx 
\cr
\end{align*}
\begin{align}
& \quad \le  c \bigg(\intl_{B(4r) }   |\bV- \bV_{ B(4r)}|^{\frac{3p}{2p-3}} dx\bigg)^{ \frac{2p-3}{3p}}
\bigg(\intl_{B(4r) \setminus B(r)}  |u|^3 \psi^3   dx \bigg) ^{ \frac{1}{3}}\bigg(\intl_{B(4r) \setminus B(r)}  |\nabla u|^p    dx \bigg) ^{ \frac{1}{p}}
\cr
& \quad+   c r^{ -1}  \bigg(\intl_{B(4r) }   |\bV- \bV_{ B(4r)}|^{\frac{3p}{2p-3}} dx\bigg)^{ \frac{2p-3}{3p}}
\bigg(\intl_{B(4r) \setminus B(r)}  |u|^3 \psi^3   dx \bigg) ^{ \frac{1}{3}}\bigg(\intl_{B(4r) \setminus B(r)}  | u|^p    dx \bigg) ^{ \frac{1}{p}}
\cr
&\quad \le  c r^{\frac23}
\bigg(\intl_{B(4r) \setminus B(r)}  |u|^3 \psi^3   dx \bigg) ^{ \frac{1}{3}}\bigg(\intl_{B(4r) \setminus B(r)}  |\nabla u|^p    dx \bigg) ^{ \frac{1}{p}}
\cr
&\qquad+ c r^{ -\frac13} \bigg(\intl_{B(4r) \setminus B(r)}  |u|^3 \psi^3   dx \bigg) ^{ \frac{1}{3}}\bigg(\intl_{B(4r) \setminus B(r)}  | u|^p    dx \bigg) ^{ \frac{1}{p}} .
 \label{2.9}
\end{align}

 Let us define  $ \widetilde{u}_{ B(4r) \setminus B(r)} = \frac{1}{ \intl \psi  dx}  \intl_{B(4r) \setminus B(r)} u \psi  dx$. Recalling that 
$ u = \nabla \cdot (\bV- \bV_{ B(2r)})$, using integration by parts, H\"older's inequality, together with  \eqref{1.2} we get
\begin{align}
|\widetilde{u}_{ B(4r) \setminus B(r)}|  &\le   \frac{1}{ \intl \psi   dx} \bigg| \intl_{B(4r) \setminus B(r)}  (\bV- \bV_{ B(4r)})\cdot \nabla \psi  dx\bigg|
\cr
& = c r^{ -1} \intmw_{B(4r)} |\bV- \bV_{ B(4r)}| dx \le  c r^{ -1} \bigg(\intmw_{B(4r)} |\bV- \bV_{ B(4r)}|^{\frac{3p}{2p-3}} dx\bigg)^{ \frac{2p-3}{3p}} 
\cr
& \le c r^{  \frac{9-7p}{3p}}. 
 \label{2.8}
\end{align}
By the triangular inequality we have
\begin{align*}
\bigg(\intl_{B(4r) \setminus B(r)}  | u|^p    dx \bigg) ^{ \frac{1}{p}} &\le \bigg(\intl_{B(4r) \setminus B(r)}  | u- u_{ B(4r) \setminus B(r)} |^p dx  \bigg) ^{ \frac{1}{p}} \cr
&\qquad  + 
\bigg(\intl_{B(4r) \setminus B(r)}  | u_{ B(4r) \setminus B(r)} - \widetilde{u}_{ B(4r) \setminus B(r)} |^p dx  \bigg) ^{ \frac{1}{p}} \cr
&\qquad +    \bigg(\intl_{B(4r) \setminus B(r)} | \widetilde{u}_{ B(4r) \setminus B(r)} |^p dx  \bigg) ^{ \frac{1}{p}}\cr
&=I_1+I_2 +I_3.
\end{align*}
Using the Poincar\'{e} inequality and \eqref{2.8},  we find 
\be\label{2.8a}
I_1+I_3 \le c r\bigg(\intl_{B(4r) \setminus B(r)}  |\nabla u|^p    dx \bigg) ^{ \frac{1}{p}} + cr^{ \frac{18-7p}{3p}}.
\ee
For $I_2$ we  use the H\"{o}lder inequality, and then the Poincar\'{e} inequality to estimate
\begin{align}\label{2.8b} 
I_2&=     \bigg(\intl_{B(4r) \setminus B(r)}   \bigg|  \frac{1}{ \intl \psi  dx}  \intl_{B(4r) \setminus B(r)}  (u- u_{ B(4r) \setminus B(r)})\psi dx\bigg|^p dx  \bigg) ^{ \frac{1}{p}} \cr
 &\le  c  \bigg(\intl_{B(4r) \setminus B(r)}  | u- u_{ B(4r) \setminus B(r)} |^p dx  \bigg) ^{ \frac{1}{p}} \le  c r\bigg(\intl_{B(4r) \setminus B(r)}  |\nabla u|^p    dx \bigg) ^{ \frac{1}{p}}. 
\end{align} 
Combining \eqref{2.8a} and \eqref{2.8b}, we get
\be\label{2.8c} 
\bigg(\intl_{B(4r) \setminus B(r)}  | u|^p    dx \bigg) ^{ \frac{1}{p}}  \le cr^{ \frac{18-7p}{3p}}+ c r\bigg(\intl_{B(4r) \setminus B(r)}  |\nabla u|^p    dx \bigg) ^{ \frac{1}{p}}.
\ee
Inserting \eqref{2.8c}  into the last term of  \eqref{2.9}  and the dividing result by $ \bigg(\intl_{B(4r) \setminus B(r)}  |u|^3 \psi^3   dx \bigg) ^{ \frac{1}{3}}$, we find 
\begin{align*}
 r^{-1}  \intl_{B(4r) \setminus B(r)} |u|^3 \psi^3   dx   
 \le cr^{-\frac13} \bigg(\intl_{B(4r) \setminus B(r)}  |\nabla u|^p    dx \bigg) ^{ \frac{1}{p}} + cr^{ \frac{18-11p}{3p}}. 
\end{align*}
Thus,    observing \eqref{2.4} and $p>9/5$,  we obtain the claim  \eqref{2.6}.   \\

Let $ 1 < r< +\infty$ be arbitrarily chosen.   By  $ \zeta \in  C^\infty(\R^{n} )$ we denote  a cut off function, which is radially non-increasing   with $ \zeta =1$ on $B(2r)$ and $\zeta=0$ on $ \Bbb R^3\setminus B(3r)$ such that  
$ |\nabla \zeta | \le c r^{ -1}$.  We multiply  \eqref{ns-steady} by 
$ u \zeta $ integrate over $ B(3r)$   and 
apply integration by parts. This yields
\begin{align}
 \intl_{B(3r) } |\nabla u|^p \zeta ^2 dx &= \intl_{B(3r)}  |\nabla u|^{p-2} \nabla  \zeta ^2 \cdot  \nabla  u  \cdot u dx \cr
 &\qquad + \frac{1}{2}\intl_{B(3r) } |u|^2 u \cdot \nabla \zeta  + \intl_{B(3r) } (\pi -\pi_{ B(3r)}) u \cdot \nabla \zeta  dx 
 \cr
 &\le   c \int_{B(3r)\setminus B(r)} |\nabla u|^p dx + c r^{-p} \int_{B(3r)\setminus B(r)} |u|^p dx \cr
 &\qquad +  c r^{ -1} \intl_{B(3r) \setminus B(2r)} |u|^3  dx
+ c  r^{ -1} \intl_{B(3r) \setminus B(2r) } |\pi - \pi_{ B(3r) }| |u|  dx\cr
&= I+ II+ III+IV. 
 \label{2.7}
\end{align}
Using    \eqref{2.5},  we immediately get 
\[
I = o(1)\quad  \text{as}\quad r \rightarrow +\infty. 
\]
From  \eqref{2.8c} and \eqref{2.4}  it follows that 
\begin{align}
II &= c\left\{  r^{-1} \bigg(\int_{B(3r)\setminus B(r)} |u|^p dx\bigg)^{\frac1p} \right\}^p\cr
 &\le c r^{ \frac{18-10p}{3}} +  c \int_{B(3r)\setminus B(r)} |\nabla u|^p dx=o(1) \quad  \text{as}\quad r \rightarrow +\infty. 
\end{align}

From  \eqref{2.6}  we also find $ III= o(1)$ as $ r \rightarrow +\infty$. Finally, applying H\"older's inequality and using  \eqref{2.6}, we get 
\begin{align}
\label{last}
IV& \le c \bigg( r^{ -1} \intl_{B(3r)} |\pi - \pi_{ B(3r)}|^{ \frac{3}{2}} dx  \bigg)^{ \frac{2}{3}} \bigg( r^{-1}\intl_{B(3r)\setminus B(r) } |u|^{3} dx  \bigg)^{ \frac{1}{3}} \cr
&=   c \bigg( r^{ -1} \intl_{B(3r)} |\pi - \pi_{ B(3r)}|^{ \frac{3}{2}} dx  \bigg)^{ \frac{2}{3}}     o(1)
\end{align}
as $r\to +\infty$.
Using the estimate  \eqref{2.21b}  with  $B( 3r)$ in place of $B( \overline{R})$, we obtain 
\begin{align*}
   r^{ -1} \intl_{B(3r)} |\pi - \pi_{ B(3r)}|^{ \frac{3}{2}} dx  &\le cr^{ \frac{9-5p}{2p}} \bigg( \intl_{B(3r)}  |\nabla u|^p dx \bigg)^{ \frac{3(p-1)}{2p}} + 
   c r^{ -1}  \intl_{B(3r)} |u|^3 dx.
\end{align*}
By virtue of  \eqref{2.4} and  \eqref{2.5} the right-hand side of the above inequality is bounded for  $ r \ge 1 $.  Therefore,  \eqref{last}  shows that $ IV=o(1)$ as $ r 
\rightarrow +\infty$. Inserting  the above estimates of $ I, II, III$ and $ IV$ into the right-hand side of  \eqref{2.7}, we deduce that 
\[
\intl_{B(r)} |\nabla u|^p dx = o(1)\quad   \text{as}\quad  r \rightarrow +\infty. 
\]
Accordingly, $ u \equiv const$ and by means of  \eqref{2.5} it follows $ u \equiv 0$.  \hfill \Beweisende

\hspace{0.5cm}
$$\mbox{\bf Acknowledgements}$$
Chae was partially supported by NRF grants 2016R1A2B3011647, while Wolf has been supported 
supported by NRF grants 2017R1E1A1A01074536.
The authors declare that they have no conflict of interest.
\bibliographystyle{siam}

\end{document}